\documentclass[a4paper,12pt]{amsart}
\usepackage{fullpage}
\usepackage{amsmath}
\usepackage{amsfonts}
\usepackage{amssymb}

\usepackage[OT4]{fontenc}
\usepackage[latin2]{inputenc}
\DeclareMathOperator\dc{{\it i\partial\overline{\partial}}}
\DeclareMathOperator\co{{\it\mathbb C^n}}

\DeclareMathOperator\we{\wedge}
\DeclareMathOperator\om{{\it\omega}}
\DeclareMathOperator\Om{{\it\Omega}}

\DeclareMathOperator \wph {{\it\om_{\phi}}}

\DeclareMathOperator \psex {{\it PSH(X,\om)}}

\DeclareMathOperator \tx {{\it TX}}
\DeclareMathOperator \tzx {{\it T_{z}X}}
\DeclareMathOperator \con{{\it D}}
\DeclareMathOperator \lbr {{\it\lbrace}}
\DeclareMathOperator \rbr {{\it\rbrace}}

\DeclareMathOperator{\tens}{\otimes}
\DeclareMathOperator\pa{\partial}

\def\om{\omega}

\def\Om{\Omega}
\def\The{\Theta}
\DeclareMathOperator \la {{\it\lambda}}
\DeclareMathOperator \eps {{\it\epsilon}}
\DeclareMathOperator \pzw {{\it\Phi(z,w)}}
\newtheorem{theorem}{Theorem}

\newtheorem{conjecture}[theorem]{Conjecture}
\newtheorem{lemma}[theorem]{Lemma}
\newtheorem{definition}[theorem]{Definition}
\newtheorem{corollary}[theorem]{Corollary}
\newtheorem{proposition}[theorem]{Proposition}

\newtheorem{remark}[theorem]{Remark}
\title{H\"older continuous potentials on manifolds with partially positive curvature}
\author{S\l awomir Dinew}

\begin{document}
\maketitle
\begin{abstract} It is proved that solutions of the complex Monge-Amp\`ere equation on compact K\"ahler manifolds with right hand side in $L^p,\ p>1$ are uniformly H\"older continuous under the assumption on non-negative orthogonal bisectional curvature.
\end{abstract}
\section{Introduction}
The regularity of weak solutions of the complex Monge-Amp\`ere euation in the setting of compact K\"ahler manifolds is a very timely topic (see \cite{EGZ}, \cite{K4}). A strong motivation for this branch of research are the recent results in the theory of the K\"ahler-Ricci flow (\cite{TZ},\cite{ST1},\cite{ST2},\cite{CN}, to mention a few). In particular, it is important to understand the limiting behavior of the K\"ahler-Ricci flow near finite time singularities (if the canonical divisor $K_X$\ is big) or the limit at infinity if additionally $K_X$\ is nef. In dimension $2$\ the picture is more or less clear and is very similar to the one in the minimal model program for surfaces (see \cite{CN}, \cite{TZ}). In higher dimensions, as explained in \cite{ST1} and \cite{ST2}, the successful termination of the program depends on understanding the following: 
\begin{enumerate}
\item The complex Monge-Amp\`ere equation on singular varieties,
\item The regularity of the solution of the corresponding Monge-Amp\`ere  equation in the limiting case.
\end{enumerate}
The first task was studied recently in \cite{EGZ}. The difficulty with the second one is due to the fact that in the limiting cases the underlying metric is no more K\"ahler and thus the equation is much less tractable. Even in the {\it big} form case the continuity of solutions is still an unsolved problem, except in some special cases (see \cite{DZ}).

However, even in the case of a K\"ahler metric the picture is far from being clear.  Thus in this note we restrict our attention to the K\"ahler case.

 The first natural question is obtaining H\"older regularity for the potentials. The basic result regarding this problem was obtained by Ko\l odziej in  \cite{K4}, where the following was proved:
\begin{theorem}\label{holderregkoj} Let $(X,\om)$\ be a compact $n$-dimensional K\"ahler manifold ($n\geq2$). Consider the following Monge-Amp\`ere equation:
\begin{equation}(\om+\dc \phi)^n=f\om^n,\ sup_X\phi=0,\ \phi\in\psex,
\end{equation}
where $f\geq 0$\ is a given function satisfying the additional property 
that $f\in L^p(\om^n),\ p>1$ (and the necessary normalization $\int_Xf\om^n=\int_X\om^n$). Then the solution $\phi$\ is H\"older continuous with H\"older exponent dependent on $X,\ n$\ and $p$.
\end{theorem}
In \cite{EGZ} it was shown that that under the additional assumption that $X$\ is {\it rational\ homogeneous} compact K\"ahler manifold the H\"older exponent is actually independent of $X$\ and can be taken to be $2/(2+nq+\epsilon),\ \forall\epsilon>0$\ ($q$\ is the conjugate of $p$ i.e. $1/p+1/q=1$). It is unknown whether this exponent is sharp, however one can produce an explicit example showing that the exponent cannot be bigger than $2/nq$ (see \cite{Pl}). 

There are some results, however, which suggest that on a general manifold the exponent should also depend only on $p$ and $n$.

First, the corresponding equation in a strictly pseudoconvex domain was studied in \cite{GKZ}, and in this setting the exponent is indeed independent of the domain considered (strictly speaking there is an universal lower bound $\gamma$, such that any solution is at least $\gamma$-H\"older continuous).

Secondly, in \cite{Bl1}, and \cite{Bl2}  the same equation was studied, but instead of the condition $f\geq0,\ f\in L^p(\om^n)$, the Author assumed that $f\geq 0$ and $f^{1/n}$\ is Lipschitz. Then a Lipschitz regularity for $\phi$\ was obtained (in \cite{Bl1} under the extra condition of nonnegative bisectional curvature, which was later removed in \cite{Bl2}).

Finally, in \cite{DZ} the following stability estimate was obtained (see also \cite{K2} for the original (slightly weaker) result):
\begin{theorem} Let $\phi$\ and $\psi$\ solve the equations
$$(\om+\dc\phi)^n=f\om^n,\ (\om+\dc\psi)^n=g\om^n,\ \phi,\ \psi\in\psex,$$
with $f,\ g\geq0,\ f,\ g\in L^p(\om^n),\ p>1$. If the solutions are normalized so that $sup_X(\phi-\psi)=sup_X(\psi-\phi)$\ (recall that solutions are bounded under the assumption that the right hand side is in $L^p,\ p>1$, see \cite{K1},\ \cite{K3}), then
 there exists a constant $c=c(p,\epsilon,c_0)$, where $c_0$\ is an upper bound for $||f||_p$ and $||g||_p$, such that
$$||\phi-\psi||_{\infty}\leq c||f-g||_1^{\frac 1{n+\epsilon}}, \forall \epsilon>0.$$
\end{theorem}
Note that again stability exponent is independent of the manifold.

All these results lead to the following conjecture (posed in \cite{DZ}):
\begin{conjecture}
The H\"older exponent in Ko\l odziej's theorem in \cite{K4} depends only on $n,\ p$\ but not on $X$.
\end{conjecture}

The main result in this note confirms the conjecture above for a large class of manifolds (containing the homogeneous ones, so in particular it covers the case considered in \cite{EGZ}, although we get a worse exponent). Let us state it:
\begin{theorem}
\label{mainresult} Let $(X,\om)$\ be a compact $n$-dimensional K\"ahler manifold ($n\geq2$). Assume additionally that $X$\ has non-negative orthogonal bisectional curvature. Consider the following Monge-Amp\`ere equation:
\begin{equation*}(\om+\dc \phi)^n=f\om^n,\ sup_X\phi=0,\ \phi\in\psex,
\end{equation*}
where $f\geq 0$\ is a given function satisfying  $\int_Xf\om^n=\int_X\om^n, f\in L^p(\om^n),\ p>1$. Then the solution $\phi$\ is H\"older continuous with H\"older exponent at least $\frac{1}{(nq+1+\eps)}$.
\end{theorem}
 Note that in \cite{Bl1} also an additional assumption on the bisectional curvature played a role.

 We would like to point out that the non-negativity of bisectional curvature is needed in our approach, but we don't know whether this geometric condition has any impact on the considered problem. In fact we still believe that the conjecture is true in general and, perhaps as in \cite{Bl2}, a new argument would allow to remove this technical assumption. 

 Recently in \cite{GuZ} (see also \cite{Ch}) a characterization of all compact K\"ahler manifolds with nonnegative orthogonal bisectional curvature was obtained (in terms of the universal covering spaces). The universal covering space of such a manifold should be isometrically biholomorphic to an element in one of the following two families:
\begin{equation}\tag{1} 
(\mathbb C^k,\ h_0)\times(M_1,\ h_1)\times\cdots\times(M_l,\ h_l)\times(\mathbb P^{n_1},\ \om_1)\times\cdots\times(\mathbb P^{n_r},\ \om_r),
\end{equation}
where $h_0$\ is the Euclidean metric, $h_i$\ are canonical metrics on irreducible compact Hermitian symmetric spaces $M_i$\ of rank $\geq 2$, and $\om_j$\ are K\"ahler metrics on $\mathbb P^{n_j}$\ carrying non-negative orthogonal bisectional curvature;
\begin{equation}\tag{2}
(Y,\ g_0)\times(M_1,\ h_1)\times\cdots\times(M_l,\ h_l)\times(\mathbb P^{n_1},\ \om_1)\times\cdots\times(\mathbb P^{n_r},\ \om_r),
\end{equation}
with $M_i,\ h_i,\ \om_j$\ as above and $Y$\ being either a simply connected Riemann surface with Gauss curvature negative somewhere, or a noncompact simply connected K\"ahler manifold ($dim Y\geq 2$) with non-negative orthogonal bisectional curvature and the minimum of the holomorphic sectional curvature is less than zero somewhere. In the latter case also the holomorphic sectional curvatures of $M_i$\ and $\mathbb P^{n_j}$\ are less than 
$$-min\lbr holomorphic\ sectional\ curvature\ of\ Y\rbr>0.$$
{\bf Notation.} Throughout the note $C$\ will denote different constants depending only on the relevant quantities. If there is a possibility of confusion these different $C$'s will be indexed. We will use the standard notation for (local) partial derivatives  $h_j:=\frac{\partial h}{\partial z_j},\  h_{\bar{j}}:=\frac{\partial h}{\partial \bar{z}_j}$. In order to make the note more readable for analysts we have decided to avoid the usage of Einstein summation. 
\section{Preliminaries}
Throughout the note we shall work on a fixed compact $n$-dimensional K\"ahler manifold $X$\ equipped with a fundamental K\"ahler form $\omega$\ (that is $d$-closed strictly positive globally defined form) given in local coordinates by 
$$\omega={\frac i 2}  \sum_{k,j=1}^n g_{k\overline{j}}dz^k \wedge d\overline z^j$$
 We assume that the metric is normalized so that 
$$\int_X\omega^n=1.$$ 

With $\om$\ one naturally associates the K\"ahler metric $g=\sum_{k,j=1}^n g_{k\overline{j}}dz^k\tens d\overline z^j$\ acting on the tangent bundle $\tx$ i.e. it defines a hermitian product on each tangent space $\tzx,\ z\in X$.

Let $\con$\ denote the Levi-Civita connection associated to the K\"ahler metric $g$. Let also $\The=\The(D)=D^2$\ be the associated curvature tensor. It is known (see \cite{De3}) that $i\The$\ is a section of the bundle $\mathcal C_{1,1}^{\infty}(X,\ Herm(\tx,\tx))$\ which
 analytically means that $i\The$ is (locally) a matrix with $(1,1)$-forms as entries. Explicitly, in a local basis $\frac{\pa}{\pa z_i},\ i=1,\cdots,n$\ of $\tx$\ the curvature tensor $\The$\ reads
$$i\The=i\sum_{1\leq k,l\leq n,\ 1\leq i,j\leq n}\mathcal C_{klij}dz_k\we d\bar{z}_l\otimes \frac{\pa}{\pa z_i}^{*}\otimes \frac{\pa}{\pa z_j},$$
for some coefficients $\mathcal C_{klij}\in\mathbb C$. 

 These coefficients can be computed according to the formula
$$i\The=i\bar{\pa}[(\bar{g})^{-1}\pa\bar{g}]\footnote{This formula varies in the literature (compare, for example, with \cite{Si}). The difference  comes simply from reversal of the row/column notation.},$$
where $g$\ is the local  matrix of the coefficients of the K\"ahler metric, $(\bar{g})^{-1}$\ is the inverse transposed matrix of $g$ and the computation is done componentwise (we use matrix notation it the above equation). In terms of the local coordinates the $(h,l)$-th entry is the $(1,1)$-form
\begin{equation}\label{equatforcherncoeff}
(-\sum_{p=1}^ng^{h\bar{p}}g_{l\bar{p}j\bar{k}}+\sum_{r,s,p=1}^n
g^{h\bar{s}}g^{r\bar{p}}g_{r\bar{s}\bar{k}}g_{l\bar{p}j})idz_j\we d\bar{z}_k,
\end{equation}
where, as usual, $g^{a\bar{b}}$\ is the entry of the inverse transposed matrix of $(g_{a\bar{b}})$ (the minus sign comes from the fact that we have interchanged the order of $dz_j$\ and $d\bar{z}_k$).

With the curvature tensor one associates the Chern curvature form which is a bilinear form on $\tx\otimes\tx$\ defined by
$$i\The(\zeta\otimes\eta,\ \zeta\otimes\eta):=\sum_{j,k,h,l=1}^nc_{j\bar{k}h\bar{l}}
\zeta_j\eta_h\bar{\zeta}_k\bar{\eta}_l,$$
where
\begin{equation}
c_{j\bar{k}h\bar{l}}:=-\sum_{h=1}^n\sum_{p=1}^ng_{h\bar{m}}
g^{h\bar{p}}g_{l\bar{p}j\bar{k}}+\sum_{h=1}^n\sum_{r,s,p=1}^ng_{h\bar{m}}
g^{h\bar{s}}g^{r\bar{p}}g_{r\bar{s}\bar{k}}g_{l\bar{p}j}idz_j=
-g_{l\bar{m}j\bar{k}}+\sum_{r,p=1}^n
g^{r\bar{p}}g_{r\bar{m}\bar{k}}g_{l\bar{p}j}.
\end{equation} 
Geometrically, this form arises as a contraction of the curvature tensor with the metric $g_{i\bar{j}}$.
The hermitian property of $i\The$\ gives the formulae
\begin{equation}\label{coeff}
\bar{c}_{k\bar{l}i\bar{j}}=c_{l\bar{k}j\bar{i}}.
\end{equation}
The k\"ahlerness of $g$\ gives us the commutation identities
$$g_{i\bar{j}k}=g_{k\bar{j}i},\ g_{i\bar{j}\bar{k}}=g_{i\bar{k}\bar{j}},$$
from which one easily deduces the well known fact that in the K\"ahler case
\begin{equation}
c_{j\bar{k}h\bar{l}}=R_{j\bar{k}h\bar{l}},
\end{equation}
where $R$\ denotes the bisectional curvature tensor. Therefore nonnegative bisectional curvature (in the K\"ahler case) is equivalent to Griffiths semi-positivity of the tangent bundle, while nonnegative orthogonal bisectional curvature is equivalent to partial positivity of the curvature of the tangent bundle in the sense of Demailly (see \cite{De1})\footnote{Recall that in \cite{De1} a more general situation was considered. In the non-K\"ahler case this equivalence fails.}.

Recall that 
\begin{equation*}
 PSH(X,\omega):=\lbrace \phi \in L^1(X,\omega):\dc\phi \geq -\omega,\ \phi \in\mathcal C^{\uparrow}(X) \rbrace
\end{equation*} 
where $\mathcal C^{\uparrow}(X)$\ denotes the space of upper semicontinuous functions. We call the functions that belong to $PSH(X,\omega)$ $\omega$-plurisubharmonic ($\om$-psh for short). We shall often use the handy notation $\wph:=\om+\dc\phi$.

The main topic of this note are the weak solutions of the following nonlinear problem:
\begin{equation}(\om+\dc \phi)^n=f\om^n,\ f\geq0,\ f\in L^p(\om^n),\ p>1,\ sup_X\phi=0,\ \phi\in PSH(X,\omega), 
\end{equation}
where the wedge product on the left hand side is defined using pluripotential theory (see \cite{K2}, \cite{K3}).
We have the following stability theorem for the solutions of this problem.
\begin{theorem}[\cite{EGZ}]\label{EGZ}
Let $u,\ v\in\psex$\ solve the equations
$$(\om+\dc u)^n=f\om^n,\ (\om+\dc v)^n=g\om^n, $$
where $f$\ and $g$\ are nonnegative functions with the right total integral and moreover $f,\ g\in L^p(\om^n),\ p>1$. Then there exists a constant $C$\ dependent on $X,\ p,\ s,\ \epsilon,\ ||f||_{L^p},\ ||g||_{L^p}$\ such that
$$||u-v||_{\infty}\leq C||u-v||_{L^s(\om^n)}^{\frac s{nq+s+\epsilon}}, \forall s>0,\ \epsilon>0.$$
\end{theorem}
By analyzing the proof of this theorem in \cite{EGZ} one obtains the following corollary.
\begin{corollary} In the setting as above if we additionally know that $u>v$, then the same conclusion holds if merely $g$\ is in $L^p$\ i.e. we put an assumtion {\bf only} on the Monge-Amp\`ere measure of the smaller function.
\end{corollary}
\section{Proof}
Let $\phi$\ be the given function with $(\om+\dc\phi)^n=f\om^n,\ f\in L^p(\om^n),\ p>1$. The main idea of the proof is to use an approximating technique, due to Demailly (see \cite{De1}), to produce an approximating sequence $\phi_{\epsilon}$\ which looks similarly to the convenient regularization by convolution.

If $z\in X$\ then the exponential mapping 
$$exp_z:\tzx\ni\zeta\rightarrow exp_z(\zeta)\in X,$$
is defined by $exp_z(\zeta):=\gamma(1)$, where $\gamma:[0,1]\rightarrow X$\ is the geodesic starting from $z=\gamma(0)$\ with the initial velocity $\frac {d\gamma}{dt}(0)=\zeta$. 

In $\mathbb C^n$\ endowed with the Euclidean metric $exp$\ reads
$exp_z(\zeta)=z+\zeta$. However on a general complex manifold the exponential mapping is not holomorphic. To circumwent this difficulty Demailly (in \cite{De1}) introduced the new mapping $exph$\ as follows:
\begin{definition}
$exph:\tx\ni(z,\zeta)\rightarrow exph_z(\zeta)\in X,\ \zeta\in\tzx$\ is defined by
\begin{enumerate}
\item $exph$\ is a $\mathcal C^{\infty}$\ smooth mapping;
\item $\forall z\in X,\ exph_z(0)=x\ {\rm and}\ d_{\zeta}exph(0)=Id_{\tzx}$;
\item $\forall z\in X$\ the map $\zeta\rightarrow exph_z{\zeta}$\ has a holomorphic Taylor expansion at $\zeta=0$. 
\end{enumerate}
\end{definition}
An intuitive way to see this mapping is to take the ordinary $\zeta\rightarrow exp_z(\zeta)$\ mapping, Taylor-expand it at $0$\ and erase all nonholomorphic terms in $\zeta$. By the E. Borel's theorem there exists a smooth mapping with Taylor expansion matching this erased data. Note that the mapping is not uniquely defined, hovewer any two mappings differ by a $\mathcal C^{\infty}$\ functions flat along the zero section of $\tx$\ which is naturally identified with $X$.

Of course, the $exph$\ mapping is still not holomorphic in $\zeta$ (the Taylor expansion at $0$\ may well be divergent), however it shares some propreties of holomorphic functions which will be satisfactory for our needs.

Let $\chi:\mathbb R_{+}\rightarrow \mathbb R_{+}$\footnote{In his proof Demailly made the explicit choice $\chi(t)=\frac C{(1-t)^2}exp(\frac 1{t-1})$\ for $t<1$\ with a suitable constant $C$.}\ be a cut-off function with
$\chi(t)=0$\ for $t\geq1,\ \chi(t)>0,\ {\rm for}\ t<1$, such that
\begin{equation}\label{total integral}
\int_{\co}\chi(||z||^2)d\lambda(z)=1,
\end{equation} 
 with $\lambda$\ the Lebesgue measure on $\co$.
We define $\phi_{\epsilon}$\ to be
\begin{equation}\label{phie}\phi_{\eps}(z)=\frac{1}{\epsilon^{2n}}\int_{\zeta\in\tzx}
\phi(exph_z(\zeta))\chi(\frac{|\zeta|^2}{\epsilon^2})d\la_{\om}(\zeta),\ \eps>0.
\end{equation}
Here $|\zeta|^2$\ stands for $\sum_{i,j=1}^ng_{i\bar{j}}(z)\zeta_i\bar{\zeta}_j$, and $d\la_{\om}(\zeta)$\ is the induced measure $\frac1{2^nn!}(\dc|\zeta|^2)^n$. Intuitively this corresponds to the familiar convolution with smoothing kernel (actually in the case of $\co$\ endowed with the Euclidean metric this is exactly the smoothing convolution).

Demailly introduced yet another function, namely for $w\in \mathbb C,\ |w|=\eps$ he defines $\pzw:=\phi_{\eps}(z)$. One has the following equality:
\begin{equation}\label{phizw}
\pzw=\int_{\zeta\in\tzx}
\phi(exph_z(w\zeta))\chi(|\zeta|^2)d\la_{\om}(\zeta).
\end{equation}
The reason for introducing the new variable $w$ (as far  as the Author  understands Demailly's intentions) is twofold:

First, Demailly proves that (a modification of) $\pzw$\ is subharmonic in the $w$\ variable and, since $\Phi$\ is a radial function in $w$, one obtains that (modifications of) $\phi_{\eps}$\ are increasing with $\epsilon$\ (like it is the case when convoluting ordinary plurisubharmonic functions)- a fact that is probably hard to achieve by merely working with $\phi_{\eps}$.

Secondly, Demailly's motivation was to obtain an approximation of $(1,1)$-currents and he was mainly interested in highly singular ones, i.e. those with nonvanishing Lelong numbers. By introducing the new variable Demailly made possible applications of Kiselman minimum principle so that one can kill the Lelong numbers up to certain level (see \cite{De1} for the details). Since our motivation is different (and since $\dc\phi$\ for a bounded $\phi$\ has all Lelong numbers zero) we shall at places proceed in a slightly different manner, for example putting attention to quantities that were immaterial in Demailly's approach.

Below we state the fundamental technical estimate for the (complex) Hessian of $\pzw$ (Proposition 3.8 in \cite{De1}):
\begin{proposition}\label{verytechnical} We have the following estimate:
\begin{align*}
&\dc\pzw[\theta,\eta]^2\geq\\
&\frac 1{\pi}|w|^2\int_{\co}-\chi_1(|\zeta|^2)\sum_{j,k,l,m=1}^n\frac{\pa^2\phi}
{\pa\bar{z}_l\pa z_m}(exph_z(w\zeta))(c_{j\bar{k}l\bar{m}}+\frac1{|w|^2}\delta_{jm}\delta_{kl})
\tau_j\bar{\tau}_kd\la(\zeta)\\
&-K'(|\theta||\eta|+|\eta|^2).
\end{align*}
Here $\dc$ acts on $X\times\mathbb C$ (i.e. also on the $w$\ variable) and the symbol $[\theta,\eta]^2$\ means that we compute the Levi form on the vector $(\theta,\eta)\in\tzx\times\mathbb C$\ (we follow Demailly's notation for the ease of the reader). The constants
 $c_{j\bar{k}l\bar{m}}$\ are the coefficients in the Chern curvature form (alternatively the coefficients of the bisectional curvature tensor). The vector $\tau$\ is of the form $\theta+\eta\zeta+O(|w|)$\footnote{For its definition and geometric meaning we refer to \cite{De1}, these are immaterial in our aproach.} and $\chi_1(t):=-\int_t^{\infty}\chi(u)du$. $K'$\ denotes a constant dependent merely on $X$ and finally $\delta_{jm}$\ is the Cronecker delta.
\end{proposition}
We shall not reproduce the technical details of the (computationally involved) proof- we refer to \cite{De1} (there is a similar computation in \cite{De2} made on the "ordinary" exponential mapping).
Instead we shall explain a detail that a priori seems to make such a result impossible.

Indeed, since $exph$\ is not holomorphic in general, when differentiating under the integral (one first parametrizes the tangent space $\tzx$\ so that we pull-back the computation of the "vector part" onto $\co$) one inevitably gets terms where elements like  $ \frac{\pa^2\phi}{\pa z_i\pa z_j}(exph_z(w\zeta))$\ appear! Since for a plurisubhamonic function we have information {\bf only} on the mixed complex derivatives (i.e. of type $\frac{\pa^2\phi}{\pa z_i\pa\bar{z}_j}$)\ one might feel that there is no hope of getting any estimates. 

Note hovewer that $exph_z$\ is a local diffeomorphism (when viewed as a real mapping) and its real Jacobian is invertible. So, if in a coordinate chart $exph$\ reads $exph_z=(F^{1},\cdots,F^{n})$, then the matrix
\begin{equation*}
\left(\begin{matrix}
 \frac{\pa F^{(i)}}{\pa z_j}& \frac{\pa F^{(i)}}{\pa\bar{z}_j}\\ 
\frac{\pa \bar{F}^{(i)}}{\pa z_j} & \frac{\pa\bar{F}^{(i)}}{\pa\bar{z}_j}
\end{matrix}\right) _{i,j}
\end{equation*}
is also invertible. 
In the matrix notation we have
\begin{equation}\label{matrix}
(\phi(F)_1',\cdots,\phi(F)_n',\phi(F)_{\bar{1}}',\cdots,
\phi(F)_{\bar{n}}')=(\phi_1',\cdots,\phi_n',
\phi_{\bar{1}}',\cdots,\phi_{\bar {n}}')_{F}\times\left(\begin{matrix}
 \frac{\pa F^{(i)}}{\pa z_j}& \frac{\pa F^{(i)}}{\pa\bar{z}_j}\\ 
\frac{\pa \bar{F}^{(i)}}{\pa z_j} & \frac{\pa\bar{F}^{(i)}}{\pa\bar{z}_j}
\end{matrix}\right)
.
\end{equation}
So, if 
$$\left(\begin{matrix}
 a_{ij}&b_{i\bar{j}}  \\ 
 c_{\bar{i}j}&d_{\bar{i}\bar{j}}
\end{matrix}\right):=
\left(\begin{matrix}
 \frac{\pa F^{(i)}}{\pa z_j}& \frac{\pa F^{(i)}}{\pa\bar{z}_j}\\ 
\frac{\pa \bar{F}^{(i)}}{\pa z_j} & \frac{\pa\bar{F}^{(i)}}{\pa\bar{z}_j}
\end{matrix}\right)^{-1},$$
then each of the coefficients among $a,\ b,\ c,\ {\rm and}\ d$\ is a smooth function and one can write down
$$\phi_{p}'(F)=\sum_{i=1}^n\phi(F)_ia_{ip}+\sum_{i=1}^n\phi(F)_{\bar{i}}
c_{\bar{i}p},$$
$$\phi_{\bar{p}}'(F)=\sum_{i=1}^n\phi(F)_ib_{i\bar{p}}+\sum_{i=1}^n
\phi(F)_{\bar{i}}
d_{\bar{i}\bar{p}}.$$
These formulae actually  allow one to integrate by parts. 

Thus, for example, doing this operation twice an element of type
$$\int\phi_{ps}''(F)F_{\bar{l}}^{p}{'}F_k^{s}{'}G$$
can be exchanged with a sum of elements of type $\int\phi(F)H_s$\ for some smooth function $H_s$. Now it is a matter of luck whether (by cancellation and estimation) these error terms can be made small. As Demailly's analysis shows we are lucky in this case.

Heuristically the argument above is nothing but the multilinear analogue of the well known integration-by-parts formula
$$\int u'(f(t))g(t)dt=\int u'(f(t))f'(t)(\frac{g(t)}{f'(t)})dt=
-\int u(f(t))\frac{d}{dt}(\frac{g(t)}{f'(t)})dt,$$
which is justified in the case when $f'\neq 0$.

In order to estimate further Demailly has shown a simple lemma (Lemma 4.4 in \cite{De1}). Here we prove a slightly different version, since we shall use the lemma in a bit different context (strictly speaking we are interested in the quantity $\delta(|w|)$\ in the original proof (see Theorem 4.1 in \cite{De1}), which was more or less immaterial for Demailly, whose main concern was the quantity $\la(z,|w|)$\ there).
\begin{lemma} Suppose $X$\ has non-negative orthogonal bisectional curvature  i.e.
$$\sum_{j,k,l,m=1}^n\frac{1}{2\pi}c_{j\bar{k}l\bar{m}}\tau_j
\bar{\tau}_{{k}}\xi_l\bar{\xi}_{m}\geq 0$$
for any tangent vecotrs $\tau$\ and $\xi$, satisfying $\sum_{i,j=1}^ng_{i\bar{j}}\tau_i\bar{\xi}_j=0$. At a point $z$\ we choose geodesic coordinates, so that $g_{i\bar{j}}=\delta_{ij}$, ($\delta$ is the Cronecker delta). Then there exist a constant $C$\ dependent only on the geometry of $X$ such that at $z$
$$\sum_{j,k,l,m=1}^n\frac{1}{2\pi}(c_{j\bar{k}l\bar{m}}+\frac1{|w|^2}
\delta_{jm}\delta_{kl})\tau_j
\bar{\tau}_{{k}}\xi_l\bar{\xi}_{m}+C|w|||\tau||^2||\xi||^2\geq 0$$
for all tangent vectors $\tau, \xi$\ and $|w|\in \mathbb R_{>0}$.
\end{lemma}
\begin{remark} Actually instead of merely the Chern curvature form Demailly considered the form
$$(i\The+u\tens Id_{\tx})(\tau\tens\xi,\tau\tens\xi)$$
for some semi-positive $(1,1)$-form $u$\ on $X$. In local coordinates it reads
$$\sum_{j,k,l,m=1}^n\frac{1}{2\pi}c_{j\bar{k}l\bar{m}}\tau_j
\bar{\tau}_{{k}}\xi_l\bar{\xi}_{m}+\sum_{j,k,l=1}^nu_{j\bar{k}}\tau_j
\bar{\tau}_{{k}}\xi_l\bar{\xi}_l.$$ 
This notion is close to the orthogonal bisectional curvature with conformal factor (see \cite{Si}, Chapter 5). For our application, however, it is crucial that $u$\ vanishes identically.
\end{remark}
\begin{proof} Let
$$\mu:=\sup_{||\tau||=1=||\xi||}|i\The(\tau\tens\xi,\tau\tens\xi)|,$$
(this quantity depends only on the geometry of $X$). We will show that actually one can take $C=5\mu\sqrt{\mu}$\ although we shall not use this explicit bound.

If $\mu=0$\ there is nothing to prove, so in the sequel we assume $\mu>0$. Note that, by assumption, $i\The(\tau\tens\xi,\tau\tens\xi)\geq0$\ whenever $\tau\perp\xi$. Consider two cases.\newline
{\it First case.} Suppose $|\sum_{i=1}^n\tau_i\bar{\xi}_i|>|w|\sqrt{\mu}$. Then
\begin{align*}
&i\The(\tau\tens\xi,\tau\tens\xi)+
\frac{1}{|w|^2}|\sum_i\tau_i\bar{\xi}_i|^2+
5\mu\sqrt{\mu}|w|||\tau||^2||\xi||^2\geq\\
&\geq-\mu+\mu+5\mu\sqrt{\mu}|w|||\tau||^2||\xi||^2\geq0.
\end{align*}
{\it Second case.} Let now $|\sum_{i=1}^n\tau_i\bar{\xi}_i|\leq|w|\sqrt{\mu}$. By homogenity it is enough to prove the inequality for all unit vectors $\tau$\ and $\xi$. Let $\tau=a\xi+b\vartheta$\ for some $a,\ b\in\mathbb C$ and some unit vector $\vartheta$\ perpendicular to $\xi$\footnote{In the sense that $\sum_{j=1}^n\xi_j\bar{\vartheta}_j=0$.}\ (so $|a|^2+|b|^2=1$\ and $|a|=|\sum_{i=1}^n\tau_i\bar{\xi}_i|\leq|w|\sqrt{\mu}$). Then we calculate
\begin{align*}
&i\The(\tau\tens\xi,\tau\tens\xi)+
\frac{1}{|w|^2}|\sum_i\tau_i\bar{\xi}_i|^2+
5\mu\sqrt{\mu}|w|||\tau||^2||\xi||^2\geq\\
&\geq|a|^2i\The(\xi\tens\xi,\xi\tens\xi)+
|b|^2i\The(\vartheta\tens\xi,\vartheta\tens\xi)+
2\Re(a\bar{b}i\The(\xi\tens\xi,\vartheta\tens\xi))+5\mu\sqrt{\mu}|w|
\end{align*}
The second term is nonnegative. We estimate the third one  using a polarization. So, we get
\begin{align*}
&|a|^2i\The(\xi\tens\xi,\xi\tens\xi)+
|b|^2i\The(i\The(\vartheta\tens\xi,\vartheta\tens\xi)+
2\Re(a\bar{b}i\The(\xi\tens\xi,\vartheta\tens\xi))+5\mu\sqrt{\mu}|w|
\geq\\
&\geq-|a|^2\mu-|a||b|[2|\Re(i\The(\xi\tens\xi,\vartheta\tens\xi))|
+2|\Im(i\The(\xi\tens\xi,\vartheta\tens\xi))|]+5\mu\sqrt{\mu}|w|\geq\\
&\geq-|a|^2\mu-|a||b||i\The(\frac{(\xi+\vartheta)}{\sqrt{2}}\tens\xi,
\frac{(\xi+\vartheta)}{\sqrt{2}}\tens\xi)-i\The(\frac{(\xi-\vartheta)}
{\sqrt{2}}\tens\xi,
\frac{(\xi-\vartheta)}{\sqrt{2}}\tens\xi)|-\\
&-|a||b||i\The(\frac{(\xi+i\vartheta)}{\sqrt{2}}\tens\xi,
\frac{(\xi+i\vartheta)}{\sqrt{2}}\tens\xi)-i\The(\frac{(\xi-i\vartheta)}
{\sqrt{2}}\tens\xi,
\frac{(\xi-i\vartheta)}{\sqrt{2}}\tens\xi)|+5\mu\sqrt{\mu}|w|\geq\\
&\geq -|a|^2\mu-4|a||b|\mu+5\mu\sqrt{\mu}|w|\geq -5|a|\mu+5\mu\sqrt{\mu}|w|\geq 0.
\end{align*}
(recall that $\vartheta$\ is perpendicular to $\xi$\ so all the vectors $\frac{(\xi\pm\vartheta)}{\sqrt{2}},\ \frac{(\xi\pm i\vartheta)}{\sqrt{2}}$\ are unitary).
\end{proof}
Now, following Demailly (\cite{De1}) if we apply the above inequality to each vector $\xi$\ in a basis of eigenvectors of $\dc\phi$, multiply by the corresponding nonnegative eigenvalue and take the sum, we get
\begin{align*}
&\frac{1}{2\pi}\sum_{jklm=1}^n\frac{\pa^2\phi}
{\pa\bar{z}_l\pa z_m}(exph_z(w\zeta))(c_{j\bar{k}l\bar{m}}+\\
&+\frac1{|w|^2}\delta_{jm}\delta_{kl})
\tau_j\bar{\tau}_k+\sum_{l=1}^n\frac{\pa^2\phi}
{\pa z_l\pa\bar{z}_l}(exph_z(w\zeta))C|w|||\tau||^2\geq0.
\end{align*}
(Strictly speaking one has to apply an approximation procedure to get this estimate, since a priori $\dc\phi$\ is merely a current, and we cannot apply pointwise linear algebra).

Plugging this into the Demailly's Proposition \ref{verytechnical} we infer

\begin{align*}
&\frac1{\pi}\dc\pzw[\theta,\ \eta]^2\geq\\
&\geq-2|w|^2\int_{\co}-\chi_1(|\zeta|^2)\sum_{l=1}^n\frac{\pa^2\phi}
{\pa z_l\pa\bar{z}_l}(exph_z(w\zeta))C|w|||\tau||^2d\la(\zeta)-\\
&-K'(|\theta||\eta|+|\eta|^2)\geq\\
&\geq-[2|w|^2\int_{\co}-\chi_1(|\zeta|^2)\sum_{l=1}^n\frac{\pa^2\phi}
{\pa z_l\pa\bar{z}_l}(exph_z(w\zeta))d\la(\zeta)]C|w||\theta|^2-\\
&-K''(|\theta||\eta|+|\eta|^2),
\end{align*}
where we have used $\tau=\theta+\eta\zeta+O(|w|)$. Since the element in the square brackets is uniformly bounded (roughly speaking it is comparable with the integral of the Laplacian over a ball with radius $|w|$) we obtain
\begin{equation}
\frac1{\pi}\dc\pzw[\theta,\ \eta]^2\geq-C|w||\theta|^2-K(|\theta||\eta|+|\eta|^2),
\end{equation}
for some constant $K$\ dependent only on $X$.

Let us stress that this estimate holds for plurisubharmonic function as well as for $\om$-psh ones. This is due to the fact that the estimate, as well as all the quantities appearing in it are local. Thus instead of working with $\phi$\ one can apply the construction to $\rho+\phi$\ with $\rho$\ a local potential for the K\"ahler metric $\om$. Everything will be well defined for small enough $\epsilon$. Of course we would have to take into account the difference between $\rho(z)$\ and $\rho_{\epsilon}(z)$, but, since the potential is smooth, this difference is easily estimated to be of order $\epsilon^2$\ i.e. it can be absorbed harmlessly into the estimate.

Note that the estimate shows that $\pzw+K|w|^2$\ is subharmonic in $w$.  Since it depends on $|w|$\ it must be an increasing function. Therefore $\phi_{\eps}+K\eps^2$\ is increasing with $\epsilon$.
Also, for a given bounded  $\om$-psh function $\phi$\ on a complex K\"ahler manifold with non-negative orthogonal bisectional curvature  the Demailly regularization gives us a quasi-psh function $\phi_{\eps}$ satisfying
$$\dc\phi_{\eps}\geq-(1+C\eps)\om,$$
for some constant $C$\ depending on the geometry of $X$\ (just take $\eta=0$, and $|w|=\epsilon$).

Without loss of generality we may assume $\phi\leq -1$. Then  $\phi_{\eps}\leq 0$\ for small $\epsilon$ (by construction) and $\varphi_{\eps}:=\frac{\phi_{\eps}+C_1\eps^2}{1+C\eps}$\ is a family of smooth $\om$-psh functions decreasing towards $\phi$\footnote{Similar idea was used in \cite{BK}.}.
By the corollary after Theorem \ref{EGZ} we have the estimate
\begin{equation}
sup_X|\varphi_{\eps}-\phi|\leq C||\varphi_{\eps}-\phi||_{L^1{(\om^n)}}^{1/(nq+1+\eps)}.
\end{equation}
Note hovewer that
$$||\varphi_{\eps}-\phi||_{\infty}=||\frac{\phi_{\eps}-\phi +C_1\eps^2-C\eps\phi}{1+C\eps}||_{\infty}\geq
\frac{||\phi_{\eps}-\phi||_{\infty}-C_2\eps}{1+C\eps},$$
for some $C_2$\ dependent on $X$\ and $||\phi||_{\infty}$. Similarly
$$||\varphi_{\eps}-\phi||_{L^1{(\om^n)}}\leq \frac{||\phi_{\eps}-\phi||_{L^1(\om^n)}+C_3\eps}{1+C\eps},$$
for some constant $C_3$ also dependent on $X$\ and $||\phi||_{\infty}$.

Note however that $||\phi_{\eps}-\phi||_{L^1(\om^n)}\leq C\eps^2$. Indeed, to see this, consider a double covering $\Om_i\Subset\Om_i^1,\ i=1,\cdots, M$. In each $\Om_i^1$\ we fix a K\"ahler potential $\rho_i$. For $\eps$\ small enough we have the inclusion
$$\lbr exph_{z}(\zeta)\ |\ z\in\Om_i,\ ||\zeta||\leq 1\rbr\subset\Om_i^1.$$
As we have already mentioned, instead of $\phi$\ one can make approximation with $\psi_i:=\rho_i+\phi_i$\ whenever the base point varies in $\Om_i$. Also it is enough to estimate $\int_{\Om_i}(\psi_i)_{\eps}(z)-\psi_i(z)dz$\ (recall that the error term coming from the regularization of the K\"ahler potential is of order $\eps^2$). But then an application of Jensen type formula as in Lemma 4.3 in \cite{GKZ} gives us the desired statement.

Coupling all these estimates one obtains
\begin{equation}
||\phi_{\eps}-\phi||_{\infty}\leq C\eps^{1/(nq+1+\eps)}
\end{equation}  
For some $C$\ dependent on $||\phi||_{\infty},\ X,\ ||f||_p,\ n\ {\rm and}\ p$\footnote{Actually, by Ko\l odziej's theorem $||\phi||_{\infty}$\ is controlled by the other quantities.}. Now arguing again in local chart, adding local K\"ahler potentials (so that we obtain true local plurisubharmonic functions) we may follow the proof of Lemma 4.2 in \cite{GKZ} to conclude that $\phi$\ is H\"older continuous with exponent $\frac{1}{(nq+1+\eps)}$\ (roughly speaking this follows from the fact that for plurisubharmonic functions the maximum over a small ball and the average over it are comparable). This completes the proof of H\"older continuity.

{\bf Acknowledgment.} The Author wishes to thank his advisor professor S\l awomir Ko\l odziej for many helpful discussions. He also wishes to thank Grzegorz Kapustka for enlightning discussions regarding complex geometry. This project was partially supported by Polish ministerial grant N N201 271135.

S\l awomir Dinew\\
Institute of Mathematics\\
Jagiellonian University\\
ul. \L ojasiewicza 6\\
30-348 Krak\'ow\\
Poland\\
\tt slawomir.dinew@im.uj.edu.pl
\end{document}